\documentclass[12pt]{article}
\usepackage{amsmath,amsthm}
\newenvironment{hangref}
  {\begin{list}{}{\setlength{\itemsep}{4pt}
  \setlength{\parsep}{0pt}\setlength{\leftmargin}{+\parindent}
  \setlength{\itemindent}{-\parindent}}}{\end{list}}

\begin{document}

\begin{center}

{\LARGE Urn Models and Fibonacci Series}\\[12pt]

\footnotesize
\vspace{1in}

\mbox{\large Yiyan Ni(*), Myron Hlynka(*), Percy H. Brill(*, **)}\\
* Department of Mathematics \& Statistics,\\
** Faculty of Business\\
University of Windsor, \\
Windsor, Ontario, Canada N9B 3P4\\
hlynka@uwindsor.ca
\normalsize

\end{center}



\noindent {\it Abstract:} We use probability urn models to discover some known and unknown series identities involving Fibonacci numbers. 

\bigskip


\noindent{\it Corresponding Author:} Myron Hlynka; hlynka@uwindsor.ca


\noindent\hrulefill
\newpage

\section{Introduction}
\label{sec-Introduction}

Urn models have a long history in probability dating back to Laplace. There are many types of urn models. One type of urn model consists of an urn with balls of two (or more) colors -- which we label Red and Blue. A ball is removed at random with replacement and additional balls are added according to some rule. Such urn models are referred to as Polya urn models. Common objectives are to determine the urn composition after some number of steps and the path which gives the result. Two monographs on the urn models are Johnson and Kotz (1977) and Mahmoud (2008). An article which includes a good bibliography is Mahmoud (2003). 

There are many research papers on infinite summations involving Fibonacci numbers. These include Brousseau (1968,1969), Rabinowitz (1999), Farhi (2015). 

In this paper we  state our objective using urn models, and add balls to the urn in such a manner that the Fibonacci numbers will be part of the scheme. We obtain closed expressions for some infinite sums involving Fibonacci numbers. We provide an easy tool to search for other interesting relationships. 

\section{Urn Model}

We begin by describing our specific urn model and our objective. The initial urn begins with $r$ red balls and $b$ blue balls. Assume all draws are independent and each ball is equally likely to be chosen on a given draw.  The probability of drawing a red ball on the first draw is $a_1=\dfrac{r}{r+b}$. We intend to draw a single ball from the urn with replacement until we eventually obtain a red ball. If we draw a blue ball, then we modify the urn (by adding either red or blue balls) so that the probability of drawing a red ball on step $i$ is $a_i$. If we were to only add blue balls to the urn at each step, then the probability of a blue ball on the $i$th draw could increase so quickly that it is possible that a red ball may never be drawn.
Our intention is to compute the probability that a red ball will be drawn eventually.  We will do this calculation in two different ways. One way will involve an infinite product and the other way will involve an infinite sum. Since both calculations are computing the same probability, the two results must be equal. So we will have created a mathematical identity. Though most of the identities that we present are already known, we believe that the method is new and quite appealing. Further there are many generailzaitons of urn models so this means there are many interesting further investigations using Fibonacci numbers.  \\
Let $a_i=P(\text{choose blue ball on draw }i)$.   Then\\
\begin{align*}
P(&\text{drawing a red ball within the first }n\text{ attempts}) \\
&=1-P(\text{drawing a blue ball for the first }n\text{ consecutive draws}) \\
&=1-a_1a_2a_3\dots a_n
\end{align*}
and 
\begin{align*}P(\text{drawing a red ball eventually})=1-a_1 a_2 \dots 
\end{align*}
Alternatively, we can compute the probability of drawing a red ball within the first $n$ draws by partitioning breaking the event into cases -- red on the first ball. blue then red on the first two balls, etc. So
\begin{align*}
P(&\text{drawing a red ball within the first }n\text{ attempts}) \\
&=P(R)+P(BR)+P(BBR)+ \dots +P(BB\dots BR) \\
&=(1-a_1)+a_1(1-a_2)+a_1a_2(1-a_3)+ \dots +a_1a_2 \dots a_{n-1}(1-a_n)
\end{align*}
and
\begin{align*}P(\text{drawing a red ball eventually})=(1-a_1)+\sum_{i=2}^{\infty}a_1\dots a_{i-1}(1-a_{i}) 
\end{align*}

Thus 
\begin{align}
1-a_1a_2a_3\dots a_n = (1-a_1)+a_1(1-a_2)+a_1a_2(1-a_3)+ \dots +a_1a_2 \dots a_{n-1}(1-a_n)
\end{align}
and
\begin{align}
1-a_1 a_2 \dots =(1-a_1)+\sum_{i=2}^{\infty}a_1\dots a_{i-1}(1-a_{i}) 
\end{align}

We observe that if we select $a_i$ in such a way that the cumulative product simplifies nicely, then some simple results may appear. Ratios of Fibonacci numbers lend themselves well to theses expressions. Even the single difference expression in (2) can be applied nicely to Fibonacci numbers. We also observe that, even though the expressions (1) and (2) were derived assuming that the entries were probabilities, this is not a requirement and the expressions are true for all real numbers. 
In Section 3, we give some examples using Fibonacci numbers in various ways.

\section{Fibonacci based examples}

\subsection{Example} Take $a_i=\dfrac{F_i}{F_{i+1}}$ for $i=1,2,\dots$. Note that $a_i\leq $ so the $a_i$ could represent probabilities. Then the left hand side (LHS) of (1) is \\
$1-a_1\dots a_n = 1-\dfrac{F_1}{F_2}\dfrac{F_2}{F_3}\dots \dfrac{F_n}{F_{n+1}}= 1- \dfrac{1}{F_{n+1}}$\\
As $n\rightarrow \infty$ in (2), $LHS\rightarrow LHS(2)= 1-0$. \\
Next look at the right hand side (RHS) of (2)
For our choice of $a_i$, we have
\begin{align*}
RHS=&(1-\dfrac{F_1}{F_2})+\dfrac{F_1}{F_2}(1-\dfrac{F_2}{F_3})+\dfrac{F_1}{F_2}\dfrac{F_2}{F_3}(1-\dfrac{F_3}{F_4})+\dots
\\
=&\sum_{i=1}^{\infty}\dfrac{1}{F_i}(1-\dfrac{F_i}{F_{i+1}})=\sum_{i=2}^{\infty}\dfrac{F_{i-1}}{F_i F_{i+1}}
\end{align*}
so 
\begin{align}\sum_{i=2}^{\infty}\dfrac{F_{i-1}}{F_i F_{i+1}}=1
\end{align}
In section 2, we noted that the urn model arose by adding blue balls to the urn. Although we do not need to know the contents of the urn at each step, it seems instructive to consider this. Since $a_1=F_1/F_2 =1,$ this means that the initial urn should have 1 blue ball only. Next $a_2=F_2/F_3=1/2$. So after drawing the blue ball from urn 1, the blue ball is replaced and 1 red ball is added, giving urn 2 with contents (1,1), where the first component represents the number of blue balls and the second component represents the number of red  balls.   Next $a_3=F_3/F_4=2/3$ so we must add 1 blue ball to the urn after drawing a blue ball and replacing it on draw 2. Thus urn 3  now contains (2,1) balls. 
Continuing in this way we find that Urn 4 has (3,2), urn 5 has (5,3), urn 6 has (8,5), and on on. The number of blue balls added each time is 0,1,1,2,3,5,... and the number of red balls added is 1,0,1,1,2,3,5, ... The Fibonacci nature of the pattern within the urn, and the additions to the urn, is clear. 

\subsection {Example}
Take $a_i=\dfrac{F_{i+1}}{2F_i}$ for $i=1,2,\dots$. Note that $a_i \leq 1$, so the $a_i$ could represent probailities. Then the left hand side (LHS) of (1) is \\
$1-a_1\dots a_n = 1-\dfrac{F_2}{2F_1}\dfrac{F_3}{2F_2}\dots \dfrac{F_{n+1}}{2F_{n}}= 1- \dfrac{F_{n+1}}{2^n}$\\
As $n\rightarrow \infty$, $LHS\rightarrow 1-0$. \\
Next look at the right hand side (RHS) of (2). 
For our choice of $a_i$, we have
\begin{align*}
RHS=&(1-\dfrac{F_2}{2F_1})+\dfrac{F_2}{2F_1}(1-\dfrac{F_3}{2F_2})+\dfrac{F_2}{2F_1}\dfrac{F_3}{2F_2}(1-\dfrac{F_4}{2F_3})+\dots
\\
&=\dfrac{1}{2}+\sum_{i=2}^{\infty}\dfrac{F_i}{2^{i-1}}(1-\dfrac{F_{i+2}}{2F_i})
=\dfrac{1}{2}+\sum_{i=2}^{\infty}\dfrac{F_i}{2^{i-1}}(\dfrac{2F_i-F_{i+1}}{2F_i})= \dfrac{1}{2}+\sum_{i=3}^{\infty}\dfrac{F_{i-2}}{2^i}
\end{align*}
so 
\begin{align}\sum_{i=3}^{\infty}\dfrac{F_{i-2}}{2^i}=1/2
\end{align}

\subsection {Example}
Take $a_i=\dfrac{F_{i+2}}{3F_i}$ for $i=1,2,\dots$. Note that $a_i \leq 1$. So the $a_i$ can represent probabilities. Then the left hand side (LHS) of (2) is \\
$1-a_1a_2a_3\dots = 1-\dfrac{F_3}{3F_1}\dfrac{F_4}{3F_2}\dfrac{F_5}{3F_3}\dots = 1-\lim_{n\rightarrow \infty} \dfrac{F_nF_{n+1}}{3^n}=1$\\
Next look at the right hand side (RHS) of (2). \\
 For our choice of $a_i$, we have
\begin{align*}
RHS&=(1-\dfrac{F_3}{3F_1})+\dfrac{F_3}{3F_1}(1-\dfrac{F_4}{3F_2})+\dfrac{F_3}{3F_1}\dfrac{F_4}{3F_2}(1-\dfrac{F_5}{3F_3})+\dots
\\
&=\dfrac{1}{3}+0+\sum_{i=3}^{\infty}\dfrac{F_iF_{i+1}}{3^{i-1}}(1-\dfrac{F_{i+2}}{3F_i})=\dfrac{1}{3}+\sum_{i=3}^{\infty}\dfrac{F_iF_{i+1}}{3^{i-1}}(\dfrac{3F_i-F_{i+2}}{3F_i})\\
&=\dfrac{1}{3}+\sum_{i=3}^{\infty}\dfrac{F_{i+1}F_{i-2}}{3^i}
\end{align*}
so 
\begin{align}\sum_{i=3}^{\infty}\dfrac{F_{i+1}F_{i-2}}{3^i}=2/3
\end{align}

 \subsection{Example} Take $a_i=\dfrac{F_i}{F_{i+2}}$ for $i=1,2,\dots$. Note that $a_i<1 $ so the $a_i$ could represent probabilities. Then the left hand side (LHS) of (1) is \\
$1-a_1\dots a_n = 1-\dfrac{F_1}{F_3}\dfrac{F_2}{F_4}\dots \dfrac{F_n}{F_{n+2}}= 1- \dfrac{1}{F_{n+1}F_{n+2}}$\\
As $n\rightarrow \infty$ in (2), $LHS\rightarrow  1-0$. \\
Next look at the right hand side (RHS) of (2)
For our choice of $a_i$, we have
\begin{align*}
RHS&=(1-\dfrac{F_1}{F_3})+\dfrac{F_1}{F_3}(1-\dfrac{F_2}{F_4})+\dfrac{F_1}{F_3}\dfrac{F_2}{F_4}(1-\dfrac{F_3}{F_5})+\dots
\\
&=1/2+1/3+\sum_{i=3}^{\infty}\dfrac{1}{F_iF_{i+1}}(1-\dfrac{F_i}{F_{i+2}})=5/6+\sum_{i=3}^{\infty}\dfrac{F_{i}}{F_i F_{i+2}}
\end{align*}
so 
\begin{align}\sum_{i=3}^{\infty}\dfrac{1}{F_{i} F_{i+2}}=1/6
\end{align}
   
 \subsection{Example} Take $a_i=\dfrac{F_i}{F_{i+3}}$ for $i=1,2,\dots$. Note that $a_i<1 $ so the $a_i$ could represent probabilities. Then the left hand side (LHS) of (1) is \\
$1-a_1\dots a_n = 1-\dfrac{F_1}{F_4}\dfrac{F_2}{F_5}\dfrac{F_3}{F_6}\dots \dfrac{F_1F_2F_3}{F_{n+1}F_{n+2}F_{n+3}}= 1- \dfrac{2}{F_{n+1}F_{n+2}F_{n+3}}$\\
As $n\rightarrow \infty$ in (2), $LHS\rightarrow= 1-0=1$. \\
Next look at the right hand side (RHS) of (2)
For our choice of $a_i$, we have
\begin{align*}
RHS&=(1-\dfrac{F_1}{F_4})+\dfrac{F_1}{F_4}(1-\dfrac{F_2}{F_5})+\dfrac{F_1}{F_4}\dfrac{F_2}{F_5}(1-\dfrac{F_3}{F_6})+\dfrac{F_1}{F_4}\dfrac{F_2}{F_5}\dfrac{F_3}{F_6}(1-\dfrac{F_4}{F_7})+\dots
\\
&=(1-\dfrac{1}{3})+\dfrac{1}{3}(1-\dfrac{1}{5})+\dfrac{1}{3}\dfrac{1}{5}(1-\dfrac{2}{8})+\sum_{i=4}^{\infty}\dfrac{2}{F_iF_{i+1}F_{i+2}}(1-\dfrac{F_i}{F_{i+3}})\\
&=\dfrac{59}{60}+\sum_{i=4}^{\infty}\dfrac{2}{F_iF_{i+1} F_{i+2}}\dfrac{2F_{i+1}}{F_{i+3}}= \dfrac{59}{60}+\sum_{i=4}^{\infty}\dfrac{4}{F_{i}F_{i+2}F_{i+3}}
\end{align*}
so 
\begin{align}\sum_{i=4}^{\infty}\dfrac{4}{F_{i} F_{i+2}F_{i+3}}=1/60
\end{align}    
 We also considered $a_i=\dfrac{F_i}{F_{i+k}}$ type cases for $k=4,5,\dots$. 

\subsection{Example} Take $a_i=\dfrac{F_{i+1}^2}{F_iF_{i+2}}$ for $i=1,2,\dots$. Note that $a_i $ can be greater than 1 so the $a_i$ cannot represent probabilities. But we can still use the relationships in (1) and (2). The left hand side (LHS) of (1) is \\
$1-a_1\dots a_n = 1-\dfrac{F_2^2}{F_1F_3}\dfrac{F_3^2}{F_2F_4}\dfrac{F_4^2}{F_3F_5}\dfrac{F_5^2}{F_4F_6}\dots \dfrac{F_{n+1}^2}{F_nF_{n+2}}
=1-\dfrac{F_{n+1}}{F_{n+2}}$\\
As $n\rightarrow \infty$ in (2), $LHS= 1-\phi^{-1}$, where $\phi=\dfrac{\sqrt{5}+1}{2}$. \\
Next look at the right hand side (RHS) of (2). Note that \\
$1-a_i=1-\dfrac{F_{i+1}^2}{F_iF_{i+2}}=\dfrac{F_iF_{i+2}-F_{i+1}^2}{F_iF_{i+2}}=\dfrac{(-1)^{i+1}}{F_iF_{i+2}}$ and 
$a_1a_2\dots a_{i-1}=\dfrac{F_i}{F_{i+1}}$ so
\begin{align*}
RHS&=(1-a_1)+a_1(1-a_2)+\dots +a_1\dots a_{n-1}(1-a_n)+\dots\\
&= \dfrac{1}{2}+\sum_{i=2}^{\infty} a_1\dots a_{i-1}(1-a_i)=\dfrac{1}{2}+\sum_{i=2}^{\infty} \dfrac{F_i}{F_{i+1}}\dfrac{(-1)^{i+1}}{F_iF_{i+2}}\\
&= \dfrac{1}{2}+\sum_{i=2}^{\infty}\dfrac{(-1)^{i+1}}{F_{i+1}F_{i+2}}
\end{align*}
so 
\begin{align}\sum_{i=2}^{\infty}\dfrac{(-1)^{i+1}}{F_{i+1}F_{i+2}}=1/2-\phi^{-1}.
\end{align}    
 
\subsection{Example} Take $a_i=\dfrac{F_iF_{i+2}}{F_{i+1}^2}$ for $i=1,2,\dots$. Note that $a_i $ can be greater than 1 so the $a_i$ cannot represent probabilities. This choice of $a_i$ is precisely the reciprocal of that in the previous example. The left hand side (LHS) of (1) is \\
$1-a_1\dots a_n = 1-\dfrac{F_1F_3}{F_2^2}\dfrac{F_2F_4}{F_3^2}\dfrac{F_3F_5}{F_4^2}\dfrac{F_4F_6}{F_5^2}\dots \dfrac{F_nF_{n+2}}{F_{n+1}^2}
=1-\dfrac{F_{n+2}}{F_{n+1}}$\\
As $n\rightarrow \infty$ in (2), $LHS= 1-\gamma$, where $\gamma=\dfrac{\sqrt{5}+1}{2}$. \\
Next look at the right hand side (RHS) of (2). Note that \\
$1-a_i=1-\dfrac{F_iF_{i+2}}{F_{i+1}^2}=\dfrac{F_{i+1}^2-F_iF_{i+2}}{F_{i+1}^2}=\dfrac{(-1)^i}{F_{i+1}^2}$ and 
$a_1a_2\dots a_{i-1}=\dfrac{F_{i+1}}{F_i}$ so
\begin{align*}
RHS&=(1-a_1)+a_1(1-a_2)+\dots +a_1\dots a_{n-1}(1-a_n)+\dots\\
&=-1+\sum_{i=2}^{\infty} a_1\dots a_{i-1}(1-a_i)=-1+\sum_{i=2}^{\infty} \dfrac{F_{i+1}}{F_i}\dfrac{(-1)^{i}}{F_{i+1}^2}\\
&= -1+\sum_{i=2}^{\infty}\dfrac{(-1)^{i}}{F_{i}F_{i+1}}
\end{align*}
so 
\begin{align}\sum_{i=2}^{\infty}\dfrac{(-1)^{i}}{F_{i}F_{i+1}}=2-\phi.
\end{align}    
It is curious that the final sum (9) looks very similar to (8).

\section*{Acknowledgments}

This research was funded through a grant from NSERC (Natural
Sciences and Engineering Research Council of Canada).

\section*{References}

\begin{hangref}
\item B.A. Brousseau, (1968) Fibonacci-Lucas Infinite Series Research Topic, The Fibonacci
Quarterly, 7, 211--217.
\item B. A.  Brousseau, (1969) Summation of Infinite Fibonacci Series,The Fibonacci Quarterly, 7, 143–-168.
\item B. Farhi, (2015) Summation of certain infinite Fibonacci related series, arXiv:1512.09033 [math.NT]
\item N.L. Johnson and S. Kotz, (1977) {\it Urn Models and Their Application.} John Wiley.
\item H. Mahmoud, (2003) P´olya Urn Models and Connections to Random Trees: A Review. JIRSS 2, 53--114. 
\item H. Mahmoud, (2008) {\it Polya Urn Models.} Chapman and Hall/CRC. 
\item S. Rabinowitz, (1999) Algorithmic Summation of Reciprocals of Products of Fibonacci Numbers,
37, 122--127.

\end{hangref}

\end{document}